\documentclass[12pt]{article}
\usepackage[centertags,intlimits]{amsmath}                     
\usepackage{color}\usepackage{amsfonts}                     
\usepackage{amsthm}
\usepackage{graphics}
\usepackage{graphicx}
\allowdisplaybreaks[2]
\numberwithin{equation}{section}
\newtheorem{theorem}{Theorem}[section]

\newtheorem{remark}{Remark}[section]

\providecommand{\abs}[1]{\lvert #1\rvert}

\providecommand{\ov}[1]{\overline{#1}}

\newcommand{\nc}{\newcommand}
\nc{\vb}{\mathbf{v}}
\nc{\bx}{\mathbf{x}}
\nc{\by}{\mathbf{y}}
\nc{\bz}{\mathbf{z}}
\nc{\bu}{\mathbf{u}}
\nc{\bv}{\mathbf{v}}
\nc{\ba}{\mathbf{a}}
\nc{\bs}{\mathbf{s}}
\nc{\bq}{\mathbf{q}}
\nc{\bd}{\mathbf{d}}
\nc{\bb}{\mathbf{b}}
\nc{\bc}{\mathbf{c}}
\nc{\bi}{\mathbf{i}}
\nc{\bfr}{\mathbf{r}}
\nc{\bA}{\mathbf{A}}
\nc{\R}{\mathbb R}
\nc{\bbc}{\mathbb C}
\nc{\N}{\mathbb N}
\nc{\D}{\mathbb D}
\nc{\Z}{\mathbb Z}
\nc{\F}{\mathbf F}
\nc{\bbS}{\mathbb S}
\nc{\B}{\cal B}
\nc{\br}{\bigr}
\nc{\bl}{\bigl}
\nc{\Bl}{\Bigl}
\nc{\Br}{\Bigr}
\nc{\ind}{\mathbf{1}}
\nc{\bP}{\mathbf{P}}

\title{On  large queue lengths in 
 generalised Jackson networks}
\author{Anatolii A. Puhalskii \footnote{Email: puhalski@iitp.ru}\\
 Institute for Problems in Information
Transmission}
\begin{document}



\maketitle
\sloppy
\begin{abstract}
This paper proves a large deviation principle (LDP) for the 
stationary distribution of queue lengths in a subcritical generalised Jackson
network assuming a Cramer condition on the interarrival and service times.
The deviation function is given by the quasipotential.  
\end{abstract}

\section{Introduction and summary}
\label{sec:introduction}

We are concerned with the asymptotics of large deviations of 
  stationary queue lengths
 in an ergodic generalised Jackson network.
It is common to characterise large deviations
 of the stationary distributions of Markov processes 
in terms of quasipotentials, which involve
trajectorial deviation
functions,  see
Freidlin
 and Wentzell \cite{wf2}, Shwartz  and Weiss \cite{SchWei95}.
Although the queue length process in a generalised Jackson network is not
Markovian,  trajectorial LDPs for that process  are
available in the literature, see, e.g., Puhalskii \cite{Puh07}. Thus,
a quasipotential can be defined, at least formally.
 In Theorem 2.2  in Puhalskii \cite{Puh19a} it was asserted that 
the quasipotential  governs large
 deviations of the
 stationary distribution.
Unfortunately the proof there is flawed as to be detailed
later.
The purpose  of this note is to give a
 correct proof.
Our proof builds on the exponential tightness of the stationary distribution
 established in Puhalskii \cite{PuhArx25}
 and  identifies the large deviation
(LD) limit point  as the time limit
of the idempotent distribution of the transient LD limit of the queue length 
 process. 

More specifically, the proof  implements the following strategy,
 cf. Puhalskii 
 \cite{Puh23a,PuhArx25}.
  Given a stationary stochastic process
$ X_n(t)$ with invariant measure $\mu_n$\,, we have that
 for a bounded continuous
nonnegative function $f$ on the state space,
\begin{equation}
  \label{eq:29'}
  \int f(x)^n\mu_n(dx)=\int\int f(x)^n\,\mathbf P( X_n(t)\in
  dx| X_n(0)=y)\mu_n(dy)\,.
\end{equation}
Suppose $\mu_n$ is exponentially tight and let $\Pi$ represent an LD
limit point of $\mu_n$\,.
(A review of the basics
 of
    LD convergence, which is just a different form
    of
the  LDP, and of idempotent processes, can be found in 
Puhalskii
\cite{Puh07} -- \cite{Puh25}.) 
If, as $n\to\infty$\,, the transient distribution $\mathbf P( X_n(t)\in
  dx| X_n(0)=y)$ continuously 
LD converges to certain deviability $\Pi_{y,t}(x)$ (or
  obeys the continuous LDP with deviation function $I_{y,t}(x)=-\ln\Pi_{y,t}(x)$
  for that matter)  and $\mu_n$ LD converges to
  deviability $\Pi$ (along a subsequence), then \eqref{eq:29'} implies that
  \begin{equation}
    \label{eq:33'}
    \sup_x f(x)\Pi(x)=\sup_y\sup_x f(x)\Pi_{y,t}(x)\Pi(y)\,.
  \end{equation}
If $\Pi_{y,t}(x)$ converges to  limit $\hat \Pi(x)$\,, as
$t\to\infty$\,, which does not depend  on $y$\,, then letting $t\to\infty$
in \eqref{eq:33'} and noting that $\sup_y\Pi(y)=1$
implies, modulo certain technical conditions,
that $\sup_x f(x)\Pi(x)=\sup_x f(x)\hat\Pi(x)$ so that
$\Pi(x)=\hat\Pi(x)$\,.
Thus, 
\begin{equation}
  \label{eq:20}
  \Pi(x)=\lim_{t\to\infty}\Pi_{y,t}(x)\,. 
\end{equation}
In the setup of this paper, 
 $X_n(t)=Q(nt)/n$\,, where $Q(t)$ represents the queue length vector
 at $t$\,. In order to be able to work with the conditional
 distribution on the righthand side of \eqref{eq:29'}, we  need
 to append $X_n(t)$ with a vector 
 of initial residual service times and initial excess exogenous
 arrival times. 
 Thus,   $y$ is a
$3K$-vector that  records the initial queue lengths at the stations as well as
the initial excess exogenous arrival times and residual service times. Therefore, the
trajectorial  LD
convergence in Puhalskii \cite{Puh07} which concerned ordinary renewal
processes needs to be extended to the case of  delayed renewal
arrival and service
processes. The same device as in Puhalskii \cite{Puh07}
of establishing exponential tightness  and
 identifying the LD limit  as a unique weak solution to a system of
idempotent equations proves its worth. With regard to taking a limit
as $t\to\infty $ in \eqref{eq:20}, the
case of
ordinary renewal processes was dealt with in
 Puhalskii \cite{Puh19a}. This paper extends the proof to the delayed
 case.
In a nutshell, since the network is subcritical, the solution of the idempotent
equations converges to a limit which implies the convergence of the solution's
idempotent distribution.

The paper is organised as follows. In Section
\ref{sec:an-ldp-stationary} a precise specification of the model is
given and the main result is stated.
In Section \ref{sec:traj} the LD trajectorial convergence for the
delayed case is addressed. The proof of the main result is discussed in
Section \ref{sec:ld-conv-stat}.
\section{An LDP for  stationary queue lengths}
\label{sec:an-ldp-stationary}

We consider a  queueing network with
 a homogeneous customer population which comprises $K$   
 single server stations.   
Customers arrive exogenously  and are served  in the
order of arrival, one customer at a time. Upon being served, they either 
join the queue at another station or leave the network.
Let $A_k(t)$  denote the   
cumulative number of exogenous arrivals  at station $k$ 
by time $t$\,,
let $S_k(t)$ denote the cumulative number of customers   
that are served   
at station $k$ for the first $t$ units of busy time of that
station, and let 
$R_{kl}(m)$
denote the cumulative number of customers among the   
first $m$ customers  departing station $k$ that go directly to station   
$l$. Let
$A_k=(A_k(t),\,t\in\R_+)$,
$S_k=(S_k(t),\,t\in\R_+)$, and $R_k=(R_{k}(m),\,m\in\Z_+)$,
where $R_{k}(m)=(R_{kl}(m),\,l\in\mathcal{K})$ and
 $\mathcal{K}=\{1,2,\ldots, K\}$\,. 
It is assumed that the $A_k$ and $S_k$ are possibly delayed 
nonzero  renewal processes
and
$  R_{kl}(m)=\sum_{i=1}^m  \mathbf1_{\{\zeta_{k}^{(i)}=l\}}$,
where $\{\zeta_k^{(1)},\zeta_k^{(2)},\ldots\}$ is  a sequence of
i.i.d. random variables assuming values in  $\mathcal{K}\cup\{0\}$\,,
$ \mathbf1_\Gamma$ standing for the indicator function of set $\Gamma$\,.
  
Let $Q_{k}(t)$ represent the number of customers present at station $k$ at time   
$t$\,. We denote $Q(t)=(Q_k(t),\,k\in\mathcal{K})$ and
$Q=(Q(t),\,t\in\R_+)$, and refer to $Q$ as to 
  the queue length process.  The random entities $A_k$\,, $S_l$\,, $R_m$ and $Q(0)$ are
assumed
to be defined on a common probability
space $(\Omega,\mathcal{F},\mathbf{P})$ and
be mutually independent,  where 
$k,l,m\in\mathcal{K}$\,.
All  the stochastic processes are assumed  to have
piecewise constant right--continuous with
left--hand limits trajectories. Accordingly, they are considered as random
elements of the associated Skorohod spaces.

We  denote $p_{kl}=\mathbf P(\zeta_k^{(1)}=l)$ and let $P=(p_{kl})_{k,l=1}^K$\,.
The matrix $P$ is assumed to be of spectral radius less than unity so
that every arriving customer eventually leaves.
For $k\in\mathcal{K}$   and $t\in\R_+$, 
 the following equations hold:
\begin{equation}
  \label{eq:2}
Q_k(t)=Q_k(0)+A_k(t)+\sum_{l=1}^KR_{lk}\bl(D_{l}(t)\br)-D_k(t),
\end{equation}
where 
\begin{equation}
  \label{eq:1}
D_k(t)=S_k\bl(B_k(t)\br)
\end{equation}
represents  the number of  departures  from station $k$ by time $t$ and
\begin{equation}
  \label{eq:5}
  B_k(t)=\int_0^t \mathbf1_{\{Q_k(u)>0\}}\,du
\end{equation}
represents the cumulative busy time of station $k$ by time $t$\,.
For  given realisations of
 $Q_k(0)$, $A_k$,  $S_k$, and $R_k$\,,
there exist unique  $Q_k=(Q_k(t),\,t\in\R_+)$, 
$D_k=(D_k(t),\,t\in\R_+)$ and $B_k=(B_k(t),\,t\in\R_+)$\,,
that satisfy 
 \eqref{eq:2}, \eqref{eq:1} and \eqref{eq:5},
where $k\in\mathcal{K}$\,, see, e.g.,
Chen and Mandelbaum \cite{MR1106809}.

Let nonnegative random variables $\xi_k$ and
$\eta_k$ represent generic times between exogenous arrivals  and
service times at station $k$, respectively. We assume that
$0<\mathbf E\xi_k<\infty$ and $0<\mathbf E\eta_k<\infty$\,.
Let  $\lambda_k=1/\mathbf E\xi_k$ and $\mu_k=1/\mathbf E\eta_k$\,.
We  assume that the network  is subcritical, i.e.
\[
\mu>(I-P^T)^{-1}\lambda
\]
where 
$\mu=(\mu_1,\ldots,\mu_K)^T$ and
$\lambda=(\lambda_1,\ldots,\lambda_K)^T$\,.
(Vector inequalities are understood  entrywise.)
In addition, it is assumed that the generic interarrival times are
unbounded and spreadout.
 More specifically, it is assumed that
 \begin{enumerate}
\item $\mathbf P(\xi_k>u)>0$\,, for all
$k\in\mathcal{K}$ and  $u>0$\,,
\item for each $k\in\mathcal{K}$\,,
there exist nonnegative function $f_k(u)$ on $\R_+$ with
  $\int_0^\infty f_k(u)\,du>0$ and $m_k\in\N$ 
such that $\mathbf P(\sum_{i=1}^{m_k}\xi_{k,i}\in[x,y])\ge
\int_x^yf_k(u)\,du$\,, provided  $0\le x\le y$\,, where 
$\xi_{k,1},\ldots,\xi_{k,m_k}$ are i.i.d. and are distributed as $\xi_k$\,.
 \end{enumerate}
Let $U(t)=(U_1(t)\,, \ldots,U_K(t))$ represent the $K$--vector of excess times of
the exogenous arrival processes at $t$\,. To put it another way, 
$U_k(t)$ is the
length of time that elapses
 between $t$  and the next exogenous arrival after $t$ at station
$k$\,.
 Let $W(t)=(W_1(t)\,, \ldots,W_K(t))$ represent the
$K$--vector of residual service times  at the
stations  at $t$\,, i.e. $W_k(t)$ is the length of time until the first
service completion after $t$\,,
with the understanding that $W_k(t)=0$ if station $k$ is idle at $t$\,.
Under the
introduced hypotheses, 
   $X(t)=(Q(t),U(t),W(t))$ is a Harris recurrent homogeneous 
Markov process, so, it converges in 
distribution to random vector $\hat X=(\hat Q,\hat U,\hat W)$\,, 
as $t\to\infty$\,, no matter the initial condition,
 see  Dai \cite{dai95}, Asmussen \cite{Asm03}.

We now recall the trajectorial LDP for the queue length process in
Puhalskii \cite{Puh07}. Suppose, in addition, that,
for all
$k\in\mathcal{K}$\,,  $A_k$ and $S_k$ are
ordinary (i.e. nondelayed) renewal processes, that
  $\mathbf{P}(\xi_k=0)=\mathbf{P}(\eta_k=0)=0,$
$\mathbf{E}\exp(\vartheta \xi_k)<\infty$\,,
$\mathbf{E}\exp(\vartheta \eta_k)<\infty$  for some 
$\vartheta>0$, and the CDFs of 
$\xi_k$ and $\eta_k$ are right differentiable at $0$ with 
 positive derivatives. 
Let
$\beta_k=\sup\{\vartheta\in\R_+:\,\mathbf{E}\exp(\vartheta
\xi_k)<\infty\}$ and
$\gamma_k=\sup\{\vartheta\in\R_+:\,\mathbf{E}\exp(\vartheta
\eta_k)<\infty\}$\,. 
Let also $\pi(u)=u\ln u-u+1$ if $u>0$\,, $\pi(0)=1$\,,
$\pi(\infty)=\infty$\,,  $0/0=0$\,, and $\infty\cdot 0=0$\,. 
Let $\mathbb{S}_+^{K\times K}$ represent the set of
row--substochastic $K\times K$--matrices with nonnegative entries and let
 $I$ represent the $K\times K$--identity matrix.
Given vectors $ a=( a_1,\ldots, a_K)^T\in\R_+^K$
and 
$ d=( d_1,\ldots, d_K)^T\in\R_+^K$,  matrix $  r \in\mathbb{S}_+^{K\times K}$ with
rows $  r _k=(r_{kl})$ and $J\subset\mathcal{K}$,
we define 
    \begin{align*}
 \psi^A_k( a_k)&=\sup_{\vartheta<\beta_k}\bl(\vartheta-  a_k\ln \mathbf{E}
\exp(\vartheta \xi_k)\br)\,,\\
\psi^S_k( d_k)&=\sup_{\vartheta<\gamma_k}\bl(\vartheta-  d_k\ln \mathbf{E}
\exp(\vartheta \eta_k)\br)\,,\\
\psi^R_k(  r _{k})&=
\sum_{l=1}^K \pi\Bl(\frac{  r _{kl}}{p_{kl}}\Br)p_{kl}+
\pi\Bl( \frac{1-\sum_{l=1}^K
    r _{kl}}{1-\sum_{l=1}^K p_{kl} }\Br)\bl(1-\sum_{l=1}^K p_{kl}\br)
\end{align*}
and
\[
     \psi_J( a, d,  r )= \sum_{k=1}^K\psi_k^A( a_k)+
\sum_{k\in J^c}\psi^S_k( d_k)
+\sum_{k\in J}\psi^S_k( d_k)
\,\mathbf1_{\{ d_k>\mu_k\}}
+\sum_{k=1}^K   \psi^R_k(  r _k)d_k\,,
\]
where  $J^c=\mathcal{K}\setminus J$\,.
 The functions $\psi_k^A( a_k)$\,,
$\psi^S_k( d_k)$ and $\psi^R_k(  r _k)$ are strictly convex,
 the functions
$\psi^S_k( d_k)
\,\mathbf1_{\{ d_k>\mu_k\}}$ are convex (note that
$\psi^S_k( \mu_k)=0$)  and
  $\psi_J( a, d,  r )\ge0$\,.
For   $y\in\R^K$, we let
\begin{equation}  \label{eq:Psi}
  L_J(y)=\inf_{\substack{( a, d,  r )\in\R_+^K\times\R_+^K\times
\mathbb{S}_+^{K\times K} 
:\\y= a+(  r ^T-I) d}}\psi_J( a, d,  r )\,.
\end{equation}
  The function $L_J(y)$ is convex by Lemma 4.2 in Puhalskii
  \cite{Puh07}. Besides, the infimum in \eqref{eq:Psi} is attained.
We also note that $L_{\mathcal{K}}(0)=0$ because one can take
$a=\lambda$\,,  $ r =P$ and $d=(I-P^T)^{-1}\lambda$ in \eqref{eq:Psi}.
On the other hand, if $J$ is a proper subset of $\mathcal{K}$\,, then 
$L_J(0)>0$\,. (Indeed, if $L_J(0)=0$\,, then $d_k=\mu_k$
and $ r _k=p_{k\cdot}$ for $k\in J^c$\,, $d_k\le\mu_k$ for $k\in J$ and
$a=\lambda$\,. Let $G$ denote the set of $k$ such that $d_k>0$\,. It
is a nonempty set because  $\lambda\not=0$\,.
If $k\in G$\,, then $ r _k=p_{k\cdot}$\,, so, for vector and matrix
restrictions to $G$\,,
$\lambda_G=(I-P^T)_{GG}d_G$\,.
Therefore,
$d_G=(I-P^T)_{GG}^{-1}\lambda_G <\mu_G$\,,
which contradicts there being $k$ with
$d_k=\mu_k$\,.)

For  $\,J\subset\mathcal{K}$\,,  
we denote 
$F_J=\{x=(x_1,\ldots,x_K)\in \R_+^K:  
x_k=0,k\in J,x_k>0,k\not\in J\}\,$\,.
Let, for
$x\in\R_+^K$ and $y\in\R^K$\,,
\begin{equation}
  \label{eq:3''}
  L(x,y)=\sum_{J\subset\mathcal{K}}
\mathbf1_{\{x\in F_J\}}L_J(y)\,.
\end{equation}
Let, given $q_0\in\R_+^K$\,,
\[
    {I}_{q_0}(q)=\int_0^\infty
L(q(t),\dot q(t))\,dt\,,
\]
provided $q\in\D(\R_+,\R_+^K)$ is  
absolutely continuous  with
$q(0)=q_0$ and $I_{q_0}(q)=\infty$\,, otherwise,
 where
$q=(q(t)\,,t\in\R_+)$ and $q(t)=(q_1(t),\ldots,q_K(t))^T$\,.
If, in addition,
\begin{equation}
  \label{eq:23}
 \lim_{n\to\infty} \mathbf
P(\abs{\frac{Q_n(0)}{n}-q_0}>\epsilon)^{1/n}=0\,, 
\end{equation}
  where
$\epsilon>0$ is  arbitrary,
 then by Theorem 2.2 in Puhalskii
\cite{Puh07} under the hypotheses  the  normalised queue length process
$\{(Q(nt)/n\,,t\in\R_+)\,,n\in\N\}$ obeys  the LDP
  in $\mathbb D(\R_+,\R^K_+)$ for rate $n$
 with deviation function $I_{q_0}(q)$\,.

For $x\in\R_+^K$\,, a quasipotential is defined by
\begin{equation}
  \label{eq:4}
  V(x)=\inf_{t\in\R_+}\inf_{\substack{q\in \D(\R_+,\R_+^K):\\
\,q(t)=x}} I_0(q)\,.
\end{equation}
(Note that in \eqref{eq:4} it can be assumed that $q(0)=0$\,.)
Our primary objective is to prove the following result.
\begin{theorem}
  \label{the:LDP}
The sequence $\{\hat Q/n\,,n\in\N\}$
obeys the LDP in $\R_+^K$ for rate $n$ with the deviation function
 $ V(x)$\,.\end{theorem}

\section{Trajectorial LD convergence}
\label{sec:traj}
The purpose of this section is to prove trajectorial LD convergence
for the delayed case. It is assumed that
\[ A_k(t)=
\ind_{\{t\ge u_k\}}(\tilde A_k(t-u_{k})+1)\]
 and that
 \[ S_k(t)=\ind_{\{t\ge v_k\}}(\tilde S_k(t-v_{k})+1)\,\]
where $\tilde A_k$ and
$\tilde S_k$ are ordinary renewal processes
and $u_k$ and $v_k$ are nonnegative numbers.
More specifically, as in Puhalskii \cite{Puh07}, we
 look at a sequence of networks indexed by $n$\,, so, some of  the pieces of
notation are appended
with subscript $n$\,.  We assume that the delays are $nu_{n,k}$ and
$nv_{n,k}$ and that time is accelerated by a factor of $n$\,.
 Thus, the model is as follows,
\begin{equation}
  \label{eq:41}
  Q_{n,k}(nt)=Q_{n,k}(0)+\ind_{\{nt\ge nu_{n,k}\}}(\tilde A_k(nt-nu_{n,k})+1)
+\sum_{l=1}^KR_{lk}\bl(D_{n,l}(nt)\br)-D_{n,k}(nt),
\end{equation}
\begin{equation}
  \label{eq:1'}
D_{n,k}(nt)=\ind_{\{B_{n,k}(nt)\ge nv_{n,k}\}}\bl(\tilde S_{k}(B_{n,k}(nt)-nv_{n,k})+1\br)\,,
\end{equation}
and \begin{equation}
  \label{eq:5a}
  B_{n,k}(nt)= \int_0^{nt} \mathbf1_{\{Q_{n,k}(u)>0\}}\,du\,,
\end{equation}
where 
$u_{n,k}\ge0$ and $v_{n,k}\ge0$\,. By convention,
$v_{n,k}>0$ if and only if $Q_{n,k}(0)>0$\,.
As a result, $B_{n,k}(nv_{n,k})=nv_{n,k}$\,.
It will be assumed that, furthermore,
$u_{n,k}\to u_k$ and $v_{n,k}\to v_k$\,.

Let
\begin{equation}
  \begin{aligned}
        \ov Q_{n,k}(t)&=\frac{Q_{n,k}(nt)}{n}\,,\quad
  \ov A_{n,k}(t)=\frac{\tilde A_{k}(nt)}{n}\,,\quad
  \ov S_{n,k}(t)=\frac{\tilde S_{k}(nt)}{n}\,,\\
  \ov R_{n,k}(t)&=\frac{R_{k}(nt)}{n}\,,\quad
  \ov D_{n,k}(t)=\frac{D_{n,k}(nt)}{n}\,,\quad
  \ov B_{n,k}(t)=\frac{B_{n,k}(nt)}{n}\,.
  \end{aligned}\label{eq:42'}
\end{equation}
Suppose that \eqref{eq:23} holds  for arbitrary $\epsilon>0$\,.
It can be proved as in Puhalskii \cite[Lemma 3.1]{Puh07} that the
stochastic processes in \eqref{eq:42'} 
are $\mathbb C$-exponentially tight of order $n$ as random elements of
$\mathbb D(\R_+,\R)$\,. The joint distributions of these processes are
$\bbc$--exponentially tight in
$\mathbb D(\R_+,\R^{6K})$\,.  The proof of the LD convergence is over if
the idempotent distribution of an LD subsequential limit of
$((\overline Q_{n,k}(t)\,, t\in\R_+)\,,k\in\mathcal{K})$  is  specified uniquely.

Let $\Pi$ represent  a deviability on  
$\mathbb D(\R_+,\R^{6K})$ that is an LD limit point of the joint
distributions  of processes in \eqref{eq:42'} for rate $n$. 
The stochastic processes 
$(\ov A_{n,k}(t)\,,t\in\R_+)$\,, $(\ov S_{n,k}(t)\,,t\in\R_+)$ and
$(\ov R_{n,k}(t)\,,t\in\R_+)$ are independent and
 LD converge to  respective independent under $\Pi$ idempotent processes
$\alpha_k=(\alpha_k(t)\,,t\in\R_+)$\,, $\sigma_k=(\sigma_k(t)\,,t\in\R_+)$ and
$\rho_{k}=(\rho_{k}(t)\,,t\in\R_+)$ which have componentwise
nondecreasing absolutely  continuous trajectories such that
 $\alpha_k(0)=0$\,, $\sigma_k(0)=0$\,,  and $\rho_{k}(0)=0$\,. The respective
 idempotent distributions of the latter idempotent processes are as
 follows, see Puhalskii \cite{Puh07},
 \begin{align}
   \label{eq:44}
   \Pi^\alpha_k(a_k)=\exp(-\int_0^\infty\psi^A_k(\dot a_k(t))\br)\,dt,\;\\
\Pi^\sigma_k(s_k)=\exp(-\int_0^\infty\psi_k^S(\dot s_k(t))\br)\,dt,\;\label{eq:44'},\\
\Pi_k^{\rho}(r_{k})=\exp(-\int_0^\infty\psi_k^R(\dot r_{k}(t))\,dt\br)\,,
\label{eq:44'''}
 \end{align}
 provided the functions $a_k=(a_k(t)\,,t\in\R_+)$\,,
$s_k=(s_k(t)\,,t\in\R_+)$ and $r_{k}=(r_{k}(t)\,,t\in\R_+)$
are nondecreasing and absolutely continuous componentwise, 
starting at $0$\,, and
$\Pi^A_k(\alpha_k=a_k)=\Pi^S_k(\sigma_k=s_k)=\Pi_k^{R}(\rho_{k}=r_{k})=0$\,, otherwise.
Let $ \chi_{k}(t)$\,,  $ \delta_{k}(t)$ and $ \beta_{k}(t)$
represent  LD subsequential limits of
  $\ov Q_{n,k}(t)$\,, $\ov D_{n,k}(t)$ and $\ov B_{n,k}(t)$,
  respectively.
These idempotent processes have
absolutely continuous trajectories and $\delta_k(t)$ and $\beta_k(t)$
are nondecreasing, with $\beta_k(t)$ growing at a rate less than or
equal to one.

The idempotent processes $(\chi_k(t)\,,t\in\R_+)$\,, $(\alpha_k(t)\,,t\in\R_+)$\,, 
$(\sigma_k(t)\,,t\in\R_+)$\,, $(\delta_k(t)\,,t\in\R_+)$\,,
$(\beta_k(t)\,,t\in\R_+)$\,, and $(\rho_{k}(t)\,,t\in\R_+)$ may be
assumed to be coordinate processes on $\bbc(\R_+,\R^{6K})$\,.
   Taking  LD subsequential limits in 
\eqref{eq:41} and \eqref{eq:1'} yields the equations
\begin{equation}
  \label{eq:43}
  \chi_k(t)=q_{0,k}+\alpha_k(t-u_{k})
+\sum_{l=1}^K\rho_{lk}\bl( \delta_{l}(t)\br)- \delta_{k}(t)
\end{equation}
and
\begin{equation}
  \label{eq:42}
  \delta_{k}(t)=\sigma_{k}(\beta_{k}(t)-v_{k})\,,
\end{equation}
where
$\alpha_k$ and $\sigma_k$ are assumed to equal zero for negative
values of the  argument and $q_{0,k}>0$ if and only if $v_k>0$\,.
In addition, passing to a subsequential LD limit in the equation
\[
  \int_0^{nt}Q_{n,k}(s)\,dB_{n,k}(s)=\int_0^{nt}Q_{n,k}(s)\,ds,  
\]
which is a consequence of \eqref{eq:5a}, implies that (see also
Puhalskii \cite{Puh07})
\begin{equation}
  \label{eq:45}
  \chi_k(t)(1-\dot \beta_k(t))=0\;\text{a.e.},
\end{equation}
where
\begin{equation}
  \label{eq:46}
  \dot \beta_k(t)\in[0,1]\;\text{a.e.}
\end{equation}
Let  $\chi(t)=(\chi_k(t)\,,k\in\mathcal{K})$\,,
$\chi=(\chi(t)\,,t\in\R_+)$\,,
$u=(u_1,\dots,u_K)$ and $v=(v_1,\ldots,v_K)$\,.
We show next that the idempotent distribution of $\chi$ is specified
uniquely and identify it.
Define, in analogy with \eqref{eq:Psi} and \eqref{eq:3''}, for
$J\subset\mathcal{K}$\,, $a=(a_1,\ldots,a_K)$\,,
$d=(d_1,\ldots,d_K)$\,, $r=(r_1,\ldots,r_K)$\,,
 $x\in\R_+$ and $y\in\R$\,,
\begin{multline}
  \label{eq:36}
  \Psi^{(J)}_{(u,v),t}( a, d,  r )= \sum_{k=1}^K\psi_k^A( a_k)\ind_{\{t> u_k\}}+
\sum_{k\in J^c}\psi^S_k( d_k)\ind_{\{t> v_k\}}
+\sum_{k\in J}\psi^S_k( d_k)
\,\mathbf1_{\{ d_k>\mu_k\}}\ind_{\{t> v_k\}}\\
+\sum_{k=1}^K   \psi^R_k(  r _k)d_k\ind_{\{t> v_k\}},
\end{multline}
\begin{equation}
  \label{eq:38}
    L^{(J)}_{(u,v),t}(y)=\inf_{\substack{( a, d,  r )\in\R_+^K\times\R_+^K \times
\mathbb{S}_+^{K\times K} 
:\\y= a+(  r ^T-I) d}}\Psi^{(J)}_{(u,v),t}( a, d,  r )
\end{equation}
and
\begin{equation}
  \label{eq:39}
  L_{(u,v),t}(x,y)=\sum_{J\subset\mathcal{K}}
\mathbf1_{\{x\in F_J\}}L^{(J)}_{(u,v),t}(y)\,.
\end{equation}
Thus, $L_{(u,v),t}(x,y)=L(x,y)$\,, for all $t$ great enough.
Let
\begin{equation}
  \label{eq:40}
  I^d_{q_0,u,v}(q)
=\int_0^\infty
L_{(u,v),t}(q(t),\dot q(t))\,dt\,,
\end{equation}
provided $q\in\mathbb D(\R_+,\R^K)$ is  
absolutely continuous  with
$q(0)=q_0$ and $I^d_{q_0,u,v}(q)=\infty$\,, otherwise.
The next theorem extends Theorem 3.1 in Puhalskii \cite{Puh07} to the
delayed case. 
\begin{theorem}
  \label{the:delayed}
The idempotent distribution of $\chi$ is given by
$\Pi(\chi=q)=\exp(-I_{q_0,u,v}^d(q))$\,.
\end{theorem}
The proof is done along the lines of the proof of Theorem 3.1 in
Puhalskii \cite{Puh07}.
Some more detail is provided next.
We replicate  the proof of   Theorem 3.1 in 
Puhalskii \cite{Puh07} up until p.120, with respective equations 
\eqref{eq:43} -- \eqref{eq:46} being used in place of 
equations (3.11) -- (3.14) in Puhalskii \cite{Puh07}.
Let $q_k(t)$ represent a sample path of $\chi_k(t)$ and let 
$a_k(t)$ represent a sample path of $\alpha_k(t-u_k)$ so that
$a_k(t)=0$\,, for $t\le u_k$\,. Similarly, let
$d_k(t)$ represent a sample path of $\delta_k(t)$\,, let $s_k(t)$ represent a
sample path of $\sigma_k(t-v_k)$\,, let $r_k(t)$ represent a sample
path of $\rho_k(t)$ and let $b_k(t)$ represent a sample
path of $\beta_k(t)$\,. Then
\begin{equation}
  \label{eq:47}
  q_k(t)=q_{0,k}+a_k(t)+\sum_{l=1}^Kr_{lk}(d_l(t))-d_k(t)\,,
\end{equation}
which is the same as equation (3.11) in Puhalskii \cite{Puh07}.
Also  \begin{equation}
\begin{aligned}
      d_k(t)&=s_k\bl(b_k(t)\br),\\
q_k(t)(1-\dot{b}_k(t))&=0\;\text{ a.e.},\\
\dot{b}_k(t)&\in[0,1]\;\text{ a.e.},
\end{aligned}\label{eq:100}
  \end{equation}
which are equations (3.12) -- (3.14) in Puhalskii \cite{Puh07}.
It is noteworthy that, as $B_{n,k}(nv_{n,k})=nv_{n,k}$\,, we have that
 $b_k(v_k)=v_k$\,.

Given a trajectory $a_k(t)$\,, we have that
\begin{multline*}
 \Pi(\cap_{t\ge u_k}\{\alpha_k(t-u_k)= a_k(t)\})=\Pi(\cap_{t\ge 0}\{\alpha_k(t)= a_k(t+u_k)\})\\
=\exp\bl(-\int_0^\infty\psi^A_k(\dot  a_k(t+u_k))\,dt\br)
=\exp(- \tilde I^A_k( a_k))\,,
\end{multline*}
where
\[
  \tilde I^A_k( a_k)=\int_{u_k}^\infty\psi^A_k(\dot{ a}_k(t))\,dt\,.
\]
Similarly,
\[
  \Pi(\cap_{t\ge v_k}\{\sigma_k(t-v_k)= s_k(t)\})=\Pi(\cap_{t\ge 
0}\{\sigma_k(t)= s_k(t+v_k)\})
=\exp(- \tilde I^S_k( s_k))\,,
\]
where, via a change of variables and the fact that $b_k(t)=t$ on 
$[0,v_k]$\,,
\begin{multline}
  \label{eq:51}
    \tilde I^S_k( s_k)=\int_{v_k}^\infty\psi^S_k(\dot
    { s}_k(t))\,dt
=\int_{0}^\infty\psi^S_k(\dot
   { s}_k(b_k(t))\dot b_k(t)\ind_{\{b_k(t)\ge v_k\}}\,dt
+\int_{b_k(\infty)}^\infty\psi^S_k(\dot
    { s}_k(t))\,dt\\=
\int_{v_k}^\infty\psi^S_k(\dot
   { s}_k(b_k(t))\dot b_k(t)\,dt
+\int_{b_k(\infty)}^\infty\psi^S_k(\dot
    { s}_k(t))\,dt\,.
\end{multline}
Analogously, 
\[
  \Pi(\cap_{t\ge0}\{\rho_k(t)= r_k(t)\})=\exp(-I^\rho_k( r_k))\,,
\]
where
\begin{multline}
  \label{eq:33}
     I^\rho_k( r_k)=\int_{0}^\infty \psi^R_k(\dot{ r}_k(t))\,dt=
\int_{v_k}^\infty \psi^R_k(\dot{ r}_k( s_k(b_k(t)))
 \dot s_k(b_k(t))\dot b_k(t)
\,dt\\+\int_{s_k(b_k(\infty))}^\infty \psi^R_k(\dot{   r}_k(t))\,dt\,.
\end{multline}
As in Puhalskii \cite{Puh07}, we let in \eqref{eq:51} and \eqref{eq:33}
$\dot s_k(t)=\mu_k$\,, for $t\ge b_k(\infty)$\,, and $\dot r_k(t)=p_k$\,,
for $t\ge s_k(b_k(\infty)$\,.

Differentiating \eqref{eq:47} yields, cf. (3.20) in Puhalskii
\cite{Puh07},
\[
  \dot q_k(t)=\dot a_k(t)+\sum_{l=1}^K\dot r_{lk}(s_l(b_l(t)))\dot
  s_l(b_l(t))\dot b_l(t)
-\dot s_k(b_k(t))\dot b_k(t)\,.
\]

Let, in analogy with the definition of $\tilde I^Q$ on p.120 in Puhalskii \cite{Puh07}
\[
  \tilde I(q)=\inf_{(a,s,r)\in \Phi(q)}(\sum_{k=1}^K\tilde I_k^A(a_k)+
\sum_{k=1}^K
\tilde  I^S_k(s_k)+\sum_{k=1}^K I^\rho_k(r_k))\,,
\]
where $\Phi(q)$ denotes the set of 
componentwise nondecreasing absolutely continuous $\R_+$--valued functions 
 $(a,s,r)$ with
$a_k(u_k)=s_k(v_k)=0$ and $r_k(0)=0$\,, for which there exist
componentwise nondecreasing absolutely continuous 
functions $b$ with  $b_k(0)=0$
    such that  \eqref{eq:47} and \eqref{eq:100} hold.
We prove that, cf. (3.23) in Puhalskii \cite{Puh07},
\begin{equation}
  \label{eq:30}
  \tilde I(q)= I^d_{q_0,u,v}(q)\,.
\end{equation}
Let as in (2.21) in Puhalskii \cite{Puh07}
\begin{align*}
  \tilde a_k(t)=\dot a_k(t)\,,&&\tilde b_k(t)=\dot  b_k(t)\,,\\
\tilde s_k(t)=\dot s_k(b_k(t))\,,
&&\tilde r_{kl}(t)=\dot r_{kl}(s_k(b_k(t))\,.
\end{align*}
Analogously to (3.24) in Puhalskii \cite{Puh07},
\begin{multline}
  \label{eq:34}
  \tilde I(q)= \inf_{\substack{(\tilde a_k(t)),(\tilde s_k(t)),\\
(\tilde r_k(t)),(\tilde b_k(t))}}\int_0^\infty \sum_{k=1}^K\bl(\psi^A_k(
 \tilde a_k(t))\ind_{\{t\ge u_k\}}
+\bl(\psi^S_k(\tilde s_k(t))
+\psi^R_k(\tilde r_k(t))\tilde s_k(t)\br)\tilde b_k(t)\ind_{\{t\ge
  v_k\}}\br)\,dt\\
=\sum_{k=1}^K\int_0^\infty\inf_{\substack{a_k,s_k,\\
r_k,b_k}}
 F_{k,t}(a_k,s_k,r_k,b_k)\,dt\,,
\end{multline}
where
\begin{equation}
  \label{eq:24}
  F_{k,t}(a_k,s_k,r_k,b_k)=
\psi^A_k(a_k)\ind_{\{t\ge u_k\}}+\bl(\psi^{S_k}(s_k)
+\psi^R_k(r_k)s_k\br)b_k\ind_{\{t\ge v_k\}}\,.
\end{equation}

Let $\mathbb{S}_+^K$ represent the set of substochastic $K$--vectors
with nonnegative entries.
The minimisation on the rightmost side of \eqref{eq:34}
is carried out over 
 $a_k\in\R_+$\,, $s_k\in\R_+$\,, 
$r_k\in \mathbb{S}_+^{ K} $ and $
b_k\in[0,1]$ subject to the constraints that, a.e.,
cf. (3.22) in Puhalskii \cite{Puh07},
 \begin{gather*}
  \dot{q}_k(t)=a_k+
\sum_{l=1}^Kr_{lk}
s_lb_l-s_kb_k,\\
q_k(t)(1-b_k)=0\,.
 \end{gather*}
Let us look at minimising
$F_{k,t}(a_k,s_k,r_k,b_k)$ over $
(a_k,s_k,r_k , b_k )$ such that, for
$k\in\mathcal{K}$\,,
\[
  y_k=a_k+\sum_{l=1}^Kr_{lk}s_lb_l-s_kb_k\,, \;\;x_k(1-b_k)=0\,,
\]
where 
$y_k\in\R$\,. With $d_k=s_kb_k$\,, 
the righthand side of \eqref{eq:24} is
\begin{equation}
  \label{eq:31}
  \sum_{k=1}^K\psi^A_k(a_k)\ind_{\{t\ge u_k\}}
+\sum_{k=1}^K\psi^S_k(\frac{d_k}{b_k})b_k\ind_{\{t\ge v_k\}}
+\sum_{k=1}^K\psi^R_k(r_k)d_k\ind_{\{t\ge v_k\}}\,.
\end{equation}
Since $\psi^{S_k}$ is a nonnegative convex function and 
$\psi^{S_k}(\mu_k)=0$\,, if $x_k>0$ and $t\ge v_k$\,,  the minimum over
$b_k$ is attained
at $b_k=1$ provided  $d_k\ge \mu_k$ and is
attained at $b_k=d_k/\mu_k$\,, provided $d_k\le\mu_k$\,.
Therefore, the minimum of the expression in \eqref{eq:31}
equals  $L^{J}_{(u,v),t}(y)$ as defined by
\eqref{eq:38}.
The equality in \eqref{eq:30} has been proved.

Let $\Pi^\chi$ represent the idempotent distribution  of
a weak solution $\chi$ to equations \eqref{eq:44}--\eqref{eq:46}.
It can be proved, as in the proof of Lemma 3.3 in Puhalskii
\cite{Puh07}, that if %
   $q$ is a piecewise linear
function, then $\Pi^\chi(q)=e^{-I_{q_0,u,v}^{d}(q)}$\,. It follows, as in the
proof of Theorem 3.1 on p.122 in Puhalskii \cite{Puh07}, that 
$\Pi^\chi(q)=e^{-I_{q_0,u,v}^{d}(q)}$ for all functions $q$\,.
Conditions (R) and (D) in Puhalskii
\cite{Puh07} needed to attend to the
technical details of the proofs are verified as in Puhalskii
\cite{Puh07}, see the proof of Theorem 2.2 on p.131 there.
This concludes the proof of Theorem \ref{the:delayed}.

As a byproduct, we obtain the next result.
\begin{theorem}
  \label{the:delayed_LDP}
Suppose that \eqref{eq:23} holds, where
$\epsilon>0$ is  arbitrary.
Then the sequence $(Q_n(nt)/n\,,
t\in\R_+)$ obeys the LDP in $\D(\R_+,\R^K)$ with deviation function 
$I^d_{q_0,u,v}(q)$\,.
\end{theorem}

\section{LD convergence of  stationary queue lengths}
\label{sec:ld-conv-stat}

In this section, Theorem  \ref{the:LDP} is proved.
We show first that the sequence 
$\hat X/n=( \hat Q/n,\hat  U/n,\hat W/n)$ is exponentially tight of order
$n$\,.
The sequence $\hat Q/n$ is exponentially tight by Puhalskii \cite{PuhArx25}.
In order to show that $\hat W/n$ is exponentially tight, we need to
prove that
\begin{equation}
  \label{eq:15}
  \lim_{u\to\infty}
\limsup_{n\to\infty}\mathbf P(\frac{\hat W}{n}>u)^{1/n}=0\,.
\end{equation}
Let $W_k(t)$ represent the residual service time of the customer in
service in server $k$ at time $t$ assuming an empty network initially
fed by an ordinary renewal process. The customer in service at $t$
enters service   with service time 
  which is independent of the past so that it has the
 generic service time distribution for station $k$\,. 
Hence,
\[  \mathbf P(\frac{W_k(t)}{n}>u)\le 
\mathbf P(\frac{\eta_{k}}{n}>u)\le e^{-\vartheta nu}\mathbf Ee^{\vartheta _k\eta_k}
\]
and
\[
  \lim_{u\to\infty}\limsup_{n\to\infty}\limsup_{t\to\infty}
\mathbf P(\frac{ W_k(t)}{n}>u)^{1/n}=0\,,
\]
proving \eqref{eq:15}. The exponential tightness of $\hat U/n$ is
proved similarly if we recall that the distribution of
 $\hat U_k$ has  density $\lambda_k\mathbf P(\xi_k>x)$\,.

Suppose now that $X(t)$ is stationary. 
Then, for $\Gamma\subset \R_+^K$\,,
\[ \mathbf
 P(\frac{Q(0)}{n}\in\Gamma)=\int_{\R_+^{3K}}\mathbf P(\frac{Q(t)}{n}
\in\Gamma|\frac{X(0)}{n}=x)\mathbf P(\frac{X(0)}{n}\in dx)\,.
\]
Hence,
\begin{equation}
  \label{eq:21}
  \mathbf P(\frac{\hat Q}{n}\in\Gamma)=
\int_{\R_+^{3K}} P_{n,x}(\frac{Q(nt)}{n}\in\Gamma)\,\mathbf P(\frac{\hat X}{n}\in dx)\,,
\end{equation}
where $P_{n,x}$ denotes the distribution on the space of trajectories
of $Q(nt)/n$
such that under $P_{n,x}$ the exogenous arrival 
processes are delayed renewal processes
with $nu$   as the vector of initial
 excess exogenous arrival times, the service time processes are
delayed renewal processes with  $nv$ as the vector of initial
residual service times and with $nq_0$ as the vector of  initial
queue lengths, where $x=(q_0,u,v)$\,.

We now check the hypotheses of  Lemma D.2 in Puhalskii
\cite{Puh_arx25} for the mixture of probabilities on the righthand
side of \eqref{eq:21}.
Suppose that $x_n=(q_{n,0},u_n,v_n)\to  x=( q_0, u, v)\in\R_+^{3K}$\,, as $n\to\infty$\,,
where $(q_{n,0},u_n,v_n)$ belongs to the support of the distribution of
$\hat X/n$\,.  
Let idempotent process $\chi_n$ represent a weak solution of the
versions of
equations
\eqref{eq:43}, \eqref{eq:42} and \eqref{eq:45} with $q_{n,0}$\,, $u_{n}$ and
$v_{n}$ as $q_{0}$\,, $u$ and $v$\,, respectively.
 Theorem \ref{the:delayed} implies that
any subsequential limit in idempotent distribution of $\chi_n$ is
the idempotent process $\chi$\,.
Hence,  the $P_{n,x_n}$ LD converge at
 rate $ n$ to the idempotent 
 distribution of $\chi$\,. We denote it by $ \Pi_{ x}$
and let $\Pi'_{ x}(y)=\Pi_{
  x}(\chi(t)=y)$\,. By \eqref{eq:40} and Theorem \ref{the:delayed},
$\Pi'_{ x}(y)=\exp\bl(-\inf_{q:\,q(t)=y}I^d_{q_0,u,v}(q)\br)$\,.
By 
\eqref{eq:36}--\eqref{eq:40}  it
is an upper semicontinuous function of $( x,y)$\,.
 The deviabilities with densities $(\Pi'_x(y)\,,y\in\R_+^K)$ are 
tight  uniformly over $x$ from compacts, i.e. given compact $K$ and
$\epsilon>0$\,, there exists compact $K'$ such that 
$\sup_{x\in K}\Pi'_x(\R_+\setminus K')<\epsilon$\,.
It is true because by \eqref{eq:43}, for $r>0$\,,
\[  \lim_{u\to\infty}\sup_{\abs{x}\le
  r}\Pi_{x}(\abs{\chi_k(t)}>u)=0\,.
\]

Assuming that  the distribution of $\hat X/n$ LD converges
along a subsequence to deviability $\hat\Pi$\,, 
we use the assertion of Lemma D.2 in Puhalskii \cite{Puh_arx25} in
order to conclude that
the righthand side
of \eqref{eq:21} LD converges to
$(\sup_{ x\in\R_+^{3K}}\Pi_{ x}(\chi(t)\in\Gamma)\hat\Pi(
 x),\Gamma\subset \R_+^K)$\,. The lefthand side LD converges to 
$(\hat\Pi(\Gamma\times
\R_+^{2K}),\Gamma\subset \R_+^K)\,$. Hence, for  bounded continuous
nonnegative function $f$ on $\R_+^K$\,,
\begin{equation}
  \label{eq:3}
  \sup_{y\in\R_+^K}f(y)\hat\Pi(\{y\}\times
\R^{2K}_+)=\sup_{y\in\R_+^K}f(y)\sup_{ x\in\R_+^{3K}}\Pi_x(\chi(t)=y)\hat\Pi(
 x)\,.
\end{equation}
We now let $t\to\infty$\,.
It can be proved as in the proof of Theorem 3.3 in Puhalskii
\cite{Puh19a} that there exists deviability $\tilde\Pi$ such that
 $\Pi_x(\chi(t) =y)\to \tilde \Pi(y)$\,, as $t\to\infty$\,, uniformly
 over $y\in\R_+^K$ and over $x$ from bounded subsets of $\R_+^{3K}$\,.
By \eqref{eq:4}, $\tilde\Pi(y)=e^{-V(y)}$\,.
By \eqref{eq:3}, on recalling that $\sup_{x\in\R_+^{3K}}\hat\Pi(x)=1$\,,
$\hat\Pi(\Gamma\times
\R_+^{2K})=\tilde\Pi(\Gamma)$\,, which concludes the proof.

Let us provide more detail on the temporal convergence. If the initial
queue length $q_0$ is a zero vector, then the function $\Pi_x(\chi(t)=y)$ is
 nondecreasing.  Indeed, by the argument of Remark 2.3 in
Puhalskii \cite{Puh19a},
$\Psi^{(\mathcal{K})}_{(u,v),t}(\lambda,\delta,P)=0$\,, where
$\delta=(I-P^T)^{-1}\lambda$\,. Hence, $L_{(u,v),t}(0,0)=0$ so that it is
possible to ''sit at the origin'' with no ''large deviation cost''.
Therefore, if $q(t)$ is an arbitrary trajectory starting at the
origin, 
then the trajectory
that stays put for some time and then replicates $q(t)$ has the same
cost.
The monotonicity implies the existence of the limit 
of $\Pi_{0,u,v}(\chi(t)=y)$\,, as $t\to\infty$\,.
As for an arbitrary initial condition $(q_0,u,v)$\,, it can be proved as
in Puhalskii \cite{Puh19a} that if $\Pi^\chi(q)\ge\kappa>0$ and 
 $(q_0,u,v)$ belongs to bounded set $B$\,, then there
exists $T>0$\,, which only depends on $\kappa$ and $B$\,, such that $q$ hits
the origin on $[0,T]$\,. (The proof relies on the network being
subcritical and is modeled on  a similar proof for the fluid limit.) It
is therefore possible to couple two different trajectories
$q(t)$\,. Assuming $t$ is great enough in order for the delays not to
matter, we can wait long enough after 
$t$ until one of the trajectories hits the origin,
keep it there at zero cost until the other trajectory gets to the
origin and glue the trajectories together when that happens. If one of
the trajectories started off at the origin, we get a coupling of an
arbitrary trajectory with a trajectory emanating from the origin so
that they have the same time limit. Details are filled in as in the
proof of Theorem 3.1 in Puhalskii \cite{Puh19a}. (Furthermore, as in
Theorem 3.1 in Puhalskii \cite{Puh19a}, for given $x$ and $y$\,,
$\Pi_x(\chi(t)=y)=\tilde\Pi(y)$\,, provided $t$ is great enough. The
time it takes to achieve that is bounded uniformly over $x$ from
bounded sets.)
\begin{remark}
 As mentioned, the proof of Theorem 2.2 in 
Puhalskii \cite{Puh19a} is flawed. 
 I misapplied Theorem 4.1 in Meyn and Down \cite{MeyDow94} 
 by assuming
 that it had to do with a standard generalised Jackson
network.
In fact, the theorem in question requires that
 the exogenous arrival processes at
the stations be obtained by splitting another  counting renewal process
so that exogenous arrivals at different stations are not independent, generally
speaking. 
\end{remark}
\def\cprime{$'$} \def\cprime{$'$} \def\cprime{$'$} \def\cprime{$'$}
  \def\cprime{$'$} \def\polhk#1{\setbox0=\hbox{#1}{\ooalign{\hidewidth
  \lower1.5ex\hbox{`}\hidewidth\crcr\unhbox0}}} \def\cprime{$'$}
  \def\cprime{$'$} \def\cprime{$'$} \def\cprime{$'$} \def\cprime{$'$}
  \def\cprime{$'$}

\end{document}